\documentclass[12pt]{article}
\usepackage{amsmath}
\usepackage{amssymb}
\usepackage{cite}
\newsymbol \blackbox 1004
\newcommand{\eh}{\hfill}\newlength{\sperr}

\newenvironment{proof}{{\settowidth{\sperr}{\bf\rm
Proof}%
\par\addvspace{0.3cm}\noindent\parbox[t]{1.3\sperr}
{\bf\rm P\eh r\eh o\eh o\eh f\eh }%
}}{\nopagebreak\mbox{}
$\blackbox$\par\addvspace{0.3cm}}

\def\nn{\nonumber}
\def\a{\alpha}
\def\b{\beta}
\def\g{\gamma}

\def\vk{\varkappa}

\def\om{\omega}

\def\t{\theta}
\def\Up{\Upsilon}
\def\vp{\varphi}

\def\ve{\varepsilon}
\def\wh{\widehat}
\def\wt{\widetilde}
\def\ov{\overline}

\def\p{\partial}

\def\BC{{\mathbb C}}
\def\BR{{\mathbb R}}
\def\BN{{\mathbb N}}
\def\clp{{\mathcal P}}

\def\cln{{\mathcal N}}
\def\clu{{\mathcal U}}
\def\cld{{\mathcal D}}
\def\cld{{\mathcal D}}
\def\cll{{\mathcal L}}
\def\im{{\rm Im\ }}

\newcommand{\E}{\mathrm{e}}
\newcommand{\I}{\mathrm{i}}
\def\mf{\mathfrak}

\newtheorem{Pa}{Paper}[section]
\newtheorem{Tm}[Pa]{{\bf Theorem}}
\newtheorem{La}[Pa]{{\bf Lemma}}
\newtheorem{Cy}[Pa]{{\bf Corollary}}
\newtheorem{Rk}[Pa]{{\bf Remark}}

\newtheorem{Dn}[Pa]{{\bf Definition}}
\newtheorem{Pn}[Pa]{{\bf Proposition}}

\title{Skew-self-adjoint Dirac systems with a rectangular matrix potential: \\
Weyl theory, direct and inverse problems}

\author{B. Fritzsche, B. Kirstein, I.Ya. Roitberg, A.L. Sakhnovich}

\date{}
\begin{document}
\maketitle

\begin{abstract} A non-classical Weyl theory is developed for
skew-self-adjoint Dirac systems with rectangular matrix potentials. The notion
of the Weyl function is introduced and 
direct and inverse problems are solved.
A Borg-Marchenko type uniqueness result and the evolution of the Weyl function
for the corresponding focusing nonlinear Schr\"odinger equation
are also derived.
\end{abstract}

{MSC(2010):}  34B20, 34L40, 37K15.

Keywords:  {\it Weyl function, Weyl theory, skew-self-adjoint Dirac system, rectangular matrix potential, direct problem, inverse problem, Borg-Marchenko type theorem, nonlinear Schr\"odinger equation.} 

\section{Introduction} \label{intro}
\setcounter{equation}{0}
Consider  a
skew-self-adjoint Dirac-type (also called a Dirac, ZS or AKNS)
system
\begin{equation}       \label{1.1}
\frac{d}{dx}y(x, z )=(\I z j+jV(x))y(x,
z ) \quad
(x \geq 0),
\end{equation}
\begin{equation}   \label{1.2}
j = \left[
\begin{array}{cc}
I_{m_1} & 0 \\ 0 & -I_{m_2}
\end{array}
\right], \hspace{1em} V= \left[\begin{array}{cc}
0&v\\v^{*}&0\end{array}\right],
 \end{equation}
where  $I_{m_k}$ is the $m_k \times
m_k$ identity
matrix and $v(x)$ is an $m_1 \times m_2$ matrix function,
which is called the potential of the system. 
Skew-self-adjoint Dirac-type systems   are very well known in
mathematics and applications, though the Weyl theory of such
systems is non-classical.
The
name ZS-AKNS comes from the fact that system (\ref{1.1}) is an
auxiliary linear system for many important nonlinear integrable
wave equations and, as such, was studied, for instance, in
\cite{AKNS, AS, CG2, FT, FKS1, GKS2, SaA1, ALS06, LAS1, ZS}. 
For  the case that $m_1\not= m_2$ systems of the form \eqref{1.1}, \eqref{1.2} are, in particular, auxiliary linear systems for
the coupled, multicomponent, and $m_1\times m_2$ matrix nonlinear
Schr\"odinger equations. We note that the above-mentioned case
$m_1\not= m_2$ was much less studied. In this respect our present paper
is a continuation of the research started in the papers \cite{FKRSp1, FKRSp2},
where self-adjoint Dirac systems with rectangular potentials were
discussed. We solve direct and inverse problems, that is, we construct Weyl functions and  
recover the $m_1 \times m_2$ potential $v$ from the Weyl function, respectively. 

Spectral and Weyl theories of auxiliary linear systems form a basis of
the  {\it Inverse Spectral Transform $($ISpT$)$  approach} \cite{Ber, BerG, KvM, SaA1, SaA2, SaA7, SaA7',
SaL1', SaL2+, SaL3} to several important initial-boundary value problems for integrable nonlinear
equations. In particular, the evolution of the Weyl functions in the form of linear fractional transformations (see \eqref{e2.9}) was first proved in \cite{SaL1', SaL2+}.
The case of the {\it focusing matrix nonlinear Schr\"odinger (fNLS) equation}
$$2v_t+\I(v_{xx}+2vv^*v)=0,$$
where $v$ are square matrices, was studied in \cite{SaA1}, and here we derive
the evolution of the Weyl functions for the fNLS with rectangular matrices.

As usual, $\BN$ stands for the set of natural numbers, $\BR$ stands for the real axis,
$\BR_+$ stands for the positive real axis,
$\BC$ stands for the complex plain, and
$\BC_+$ for the open upper
semi-plane. The notation $\BC_M$ stands for the open 
semi-plane $\{z:\, \Im (z)>M>0\}$, where $\Im (z)$ (or $\Im z$) is the imaginary part of $z$.
If $a\in \BC$, then $\Re (a)$ is its real part and $\ov a$  its complex conjugate. If $\cll$ is a set,
then    $\ov \cll$ is the closure of the set. The notation $\im$
is used for image.
An $m_2 \times m_1$ matrix $\a$ is said to be non-expansive, 
if $\a^*\a \leq I_{m_1}$ (or, equivalently, if $\a\a^* \leq I_{m_2}$).
We assume that the spaces $\BC^{m_i}$ are equipped with the $l^2$-norm
and use $\|\a\|$ to denote the corresponding operator norm of the matrix $\a$.
If the matrix function $f$ is absolutely continuous and $\sup\|f^{\prime}\|<\infty$,
we say that $f$ is a function with bounded derivative or (equivalently)
that $f$  is boundedly differentiable. The set of functions boundedly differentiable
on $[0,\, l]$ is denoted by $B^1[0,\, l]$.
We use  $I$ to denote
the identity operator and
$B({\bf H_1}, {\bf H_2})$ to denote the
class of bounded operators acting from ${\bf H_1}$ to ${\bf H_2}$. We write 
$B({\bf H_1})$ if ${\bf H_1}={\bf H_2}$.  

We set $m_1+m_2=:m$. The fundamental solution of system
\eqref{1.1} is denoted by $u(x,z)$, and this solution is normalized by the condition
\begin{align}&      \label{1.3}
u(0,z)=I_m.
\end{align}
\section{Direct problem} \label{DP}
\setcounter{equation}{0}
We consider the Dirac system \eqref{1.1} on the semi-axis $x \in [0,\, \infty)$
and assume that $v$ is measurable and locally bounded, that is, bounded
on all finite intervals.
Using a method which is similar to the method proposed in \cite{SaA1, ALS06},
where the case $m_1=m_2$ was treated,
we shall apply M\"obius transformations and matrix balls to solve the direct problem for Dirac
system. 

We next introduce a class  ${\bf P}_{\geq}(j,\BC_M)$ of nonsingular $m \times m_1$ matrix functions 
$\clp(z)$ with property-$j$, which are an immediate analog of the classical pairs
of parameter matrix functions. More specifically, the matrix functions 
$\clp(z)$ are meromorphic in $\BC_+$ and satisfy
(excluding, possibly, a discrete set of points)
the following relations
\begin{equation}\label{2.1}
\clp(z)^*\clp(z) >0, \quad \clp(z)^* j \clp(z) \geq 0 \quad \mathrm{for} \,\, z\in \BC_M= \{z:\, \Im z>M>0\}.
\end{equation}
\begin{Dn} \label{set}
The set $\cln(x,z)$ of M\"obius transformations is the set of values $($at the fixed points $x, \,z)$ 
of matrix functions
\begin{align}\label{2.2}&
\vp(x,z,\clp)=\begin{bmatrix}
0 &I_{m_2}
\end{bmatrix}u(x,z)^{-1}\clp(z)\Big(\begin{bmatrix}
I_{m_1} & 0
\end{bmatrix}u(x,z)^{-1}\clp(z)\Big)^{-1},
\end{align}
where $\clp(z)$ are nonsingular  matrix functions 
 with property-$j$, that is, 
 \begin{align}\label{2.2'}&
\clp(z)\in {\bf P}_{\geq}(j,\BC_M).
\end{align}
 \end{Dn}
 \begin{Pn} \label{PnW1} Let the Dirac system \eqref{1.1} on $[0, \, \infty)$
 be given and assume that $\|v\|$ is  bounded by $M$:
 \begin{align}\label{s1}
 \|v(x)\| \leq M  \quad \mathrm{for} \,\, x\in \BR_+.
 \end{align}
 Then the sets $\cln(x,z)$
 are well-defined  in $\BC_M$. There is a unique matrix function
 $\vp(z)$  such that
\begin{align}&      \label{2.3}
\vp(z)=\bigcap_{x<\infty}\cln(x,z).
\end{align} 
This function is analytic and non-expansive in $\BC_M$.
 \end{Pn}
\begin{proof}. The proof is similar to the proof of \cite[Proposition 2.2]{FKRSp1},
only $z$ is considered now in the semi-plane $\BC_M$ instead of the semi-plane
$\BC_+$ in \cite{FKRSp1}. We shall not, however, skip the details of this important proof,
although we will later omit the details of proofs similar to those in \cite{FKRSp1, FKRSp2}.
It is immediately apparent from \eqref{1.1} and \eqref{s1} that
\begin{align}&      \label{2.4}
\frac{d}{dx}\big(u(x,z)^*ju(x,z)\big)=-2u(x,z)^*(\Im( z) I_m-V(x))u(x,z)<0, \quad z \in \BC_M.
\end{align}
According to \eqref{1.3} and \eqref{2.4} we have
\begin{align}&      \label{2.5}
\mathfrak{A}(x,z)=\{\mathfrak{A}_{ij}(x,z)\}_{i,j=1}^2:= u(x,z)^*ju(x,z)\leq j,
\quad z \in \BC_M,
\end{align}
where $\mathfrak{A}$ is partitioned into four blocks so that $\mathfrak{A}_{ii}$
is an $m_i \times m_i$ matrix function ($i=1,2$).
Inequality \eqref{2.5} yields
\begin{align}&      \label{2.6}
\big(u(x,z)^*\big)^{-1}ju(x,z)^{-1}\geq j.
\end{align}
In view of \eqref{2.1}, we note that $\clp(z)$ maps $\BC^{m_1}$ into a maximal 
$j$-nonnegative subspace of $\BC^m$. Inequality \eqref{2.6} implies that
$u(x,z)^{-1}$ maps all maximal $j$-nonnegative subspaces into maximal 
$j$-nonnegative subspaces, that is, the image $\im \big(u(x,z)^{-1}\clp(z)\big)$
is a maximal $j$-nonnegative subspace. Clearly $\im \big(\begin{bmatrix}
I_{m_1} & 0
\end{bmatrix}^*\big)$ is a maximal $j$-nonnegative subspace too.
Thus, we get 
\begin{align}&      \label{2.7}
\det \Big(\begin{bmatrix}
I_{m_1} & 0
\end{bmatrix}u(x,z)^{-1}\clp(z)\Big)\not= 0,
\end{align}
and so $\cln$ is well-defined via \eqref{2.2}. Indeed, if \eqref{2.7} does not
hold, there is a vector $f \in \BC^{m_1}$ such that
\begin{align}&      \label{2.8}
\begin{bmatrix}
I_{m_1} & 0
\end{bmatrix}j u(x,z)^{-1}\clp(z)f=\begin{bmatrix} I_{m_1} & 0
\end{bmatrix}u(x,z)^{-1}\clp(z)f= 0, \quad f\not=0.
\end{align}
Since the left-hand side in 
\eqref{2.8} equals zero and (as discussed above)  we deal with the  maximal $j$-nonnegative
subspaces, we have
$u(x,z)^{-1}\clp(z)f \in \im \big(\begin{bmatrix}
I_{m_1} & 0
\end{bmatrix}^*\big)$. But then it follows from the second equality in \eqref{2.8}
that $f=0$, which contradicts the inequality in \eqref{2.8}.

Next, we rewrite \eqref{2.2} in the equivalent form
\begin{align}&      \label{2.9}
\begin{bmatrix}
I_{m_1} \\ \vp(x,z,\clp)
\end{bmatrix}=u(x,z)^{-1}\clp(z)\Big(\begin{bmatrix}
I_{m_1} & 0
\end{bmatrix}u(x,z)^{-1}\clp(z)\Big)^{-1}.
\end{align}
In view of \eqref{2.1}, \eqref{2.9}, and of the definition of 
$\mathfrak{A}$ in \eqref{2.5},
the formula  
\begin{align}&      \label{2.10}
\wh \vp(z) \in \cln(x,z)
\end{align}
is equivalent to
\begin{align}&      \label{2.11}
\begin{bmatrix}
I_{m_1} & \wh \vp(z)^*
\end{bmatrix}\mathfrak{A}(x,z)\begin{bmatrix}
I_{m_1} \\ \wh \vp(z)
\end{bmatrix} \geq 0.
\end{align}
In the standard way, using formula \eqref{2.4} and the equivalence of \eqref{2.10} and \eqref{2.11}, we get
\begin{align}&      \label{2.16}
\cln(x_1,z) \subset \cln(x_2,z)  \quad {\mathrm{for}} \quad x_1>x_2.
\end{align}
Moreover, \eqref{2.11} at $x=0$ means that
\begin{align}&      \label{2.17}
\cln(0,z) =\{\wh \vp(z): \, \wh \vp(z)^* \wh \vp(z)\leq I_{m_1}\}.
\end{align}
By Montel's theorem, formulas \eqref{2.16} and \eqref{2.17}
imply the existence of an analytic and non-expansive matrix function
$\vp(z)$ such that
\begin{align}&      \label{2.18}
\vp(z)\in \bigcap_{x<\infty}\cln(x,z).
\end{align} 
Indeed, because of \eqref{2.16} and \eqref{2.17} we see that the set 
of functions $\vp(x,z,\clp)$ of the form \eqref{2.2} is uniformly bounded
in $\BC_M$. So, Montel's theorem is applicable and there is an analytic
matrix function, which we denote by $\vp_{\infty}(z)$ and which is a uniform limit
of some sequence
\begin{align}&      \label{M1}
\vp_{\infty}(z)=\lim_{i \to \infty} \vp(x_i,z,\clp_i) \quad (i \in \BN, \quad x_i \uparrow, \quad \lim_{i \to \infty}x_i=\infty)
\end{align} 
on all the bounded and closed subsets of $\BC_M$. Since $x_i \uparrow$
and equalities \eqref{2.9} and \eqref{2.16} hold, it follows that the matrix functions
\begin{align}&      \nn
\clp_{ij}(z):=u(x_i,z)\begin{bmatrix}
I_{m_1} \\ \vp(x_j,z,\clp_j) 
\end{bmatrix} \quad (j \geq i)
\end{align} 
satisfy relations \eqref{2.1}. Therefore, using \eqref{M1} we see that \eqref{2.1} holds for
\begin{align}&      \nn
\clp_{i, \infty}(z):=u(x_i,z)\begin{bmatrix}
I_{m_1} \\ \vp_{\infty}(z) 
\end{bmatrix},
\end{align} 
which  implies that we can substitute $\clp=\clp_{i, \infty}$ and $x=x_i$
into \eqref{2.9} to get
\begin{align}&      \label{M2}
\vp_{\infty}(z)\in\cln(x_i,z).
\end{align} 
Since \eqref{M2} holds for all $i\in \BN$,
we see that \eqref{2.18} is true for $\vp(z)=\vp_{\infty}(z)$.

Now, let us show  that $\cln$ is a matrix ball.  It follows from
\eqref{2.4} and \eqref{2.5} that 
\begin{align}&      \label{s2}
\frac{d}{dx}\mf{A} \leq -2(\Im (z) -M)u^*u \leq 2(\Im (z) -M)\mf{A}\leq
 2(\Im (z) -M)j.
\end{align} 
Taking into account the inequalities above and the equality $\mathfrak{A}_{22}(0,z)= -I_{m_2}$, we obtain
\begin{align}&      \label{2.12}
- \mathfrak{A}_{22}(x,z) \geq \big(1+2(\Im (z) -M)x\big)I_{m_2} \quad (z\in \BC_M).
\end{align}
Note also that \eqref{2.5} implies $\mf{A}(x,z)^{-1} \geq j$ 
for $z\in \BC_M$ (see \cite{P}). Thus, we get
\begin{align}&      \label{2.13}
\big(\mathfrak{A}^{-1}\big)_{11}=\big(\mathfrak{A}_{11}-
\mf{A}_{12}\mf{A}_{22}^{-1}\mf{A}_{21}\big)^{-1} \geq I_{m_1}.
\end{align}
Since $- \mf{A}_{22}>0$, the square root $\Up=\big(- \mf{A}_{22}\big)^{1/2}$
is well-defined and we rewrite \eqref{2.11} in the form
\begin{align}&      \nn
\mathfrak{A}_{11}-
\mf{A}_{12}\mf{A}_{22}^{-1}\mf{A}_{21}-\big(\wh \vp^*\Up-\mf{A}_{12}\Up^{-1}\big)
\big(\Up \wh \vp-\Up^{-1}\mf{A}_{21}\big) \geq 0, \end{align}
where $\mf{A}_{12}=\mf{A}_{21}^*$. Equivalently, we have
\begin{align}&      \label{2.14}
\wh \vp =\rho_l \om \rho_r-\mf{A}_{22}^{-1}\mf{A}_{21}, \quad \om^*\om \leq I_{m_2},
\\ &      \label{2.15}
 \rho_l:=\Up^{-1}=\big(- \mf{A}_{22}\big)^{-1/2},
\quad \rho_r:=(\mathfrak{A}_{11}-
\mf{A}_{12}\mf{A}_{22}^{-1}\mf{A}_{21}\big)^{1/2}.
\end{align}
Here $\om$ is an $m_2 \times m_1$ non-expansive matrix function,  $\rho_l$ and $ \rho_r$ are
the left and right semi-radii of the Weyl disc.
Since \eqref{2.10} is equivalent to \eqref{2.14}, the sets
$\cln(x,z)$ (where the values of $x$ and $z$ are fixed) are matrix balls.
According to \eqref{2.12}, \eqref{2.13}, and \eqref{2.15}
the next formula holds:
\begin{align}&      \label{2.19}
\rho_l(x,z) \to 0 \quad (x \to \infty), \quad \rho_r(x,z) \leq I_{m_1}.
\end{align}
Finally, relations \eqref{2.18}, \eqref{2.14}, and \eqref{2.19}
imply \eqref{2.3}.
\end{proof}
\begin{Rk}\label{Rk3} It follows from the proof of  Proposition \ref{PnW1}
that the inequality
 \begin{align}\label{s1'}
 \|v(x)\| \leq M  \quad \mathrm{for} \quad0< x<l<\infty
 \end{align}
 implies that the embedded matrix balls $\cln(x,z)$ are well-defined for $z \in \BC_M$
 and $x$ on the interval $[0, \, l]$.
\end{Rk}
\begin{Rk}\label{Rk4}  If $\vp(z)\in \cln(l,z)$  for all $z \in \BC_M$ and, in addition,
$\vp$ is analytic in $\BC_M$, then it admits representation \eqref{2.2}, where $x=l$ and $\clp$
satisfies \eqref{2.2'}. This statement is immediately apparent, since we can choose
\[
\clp(z)=u(l,z)\begin{bmatrix}
I_{m_1} \\ \vp(z)
\end{bmatrix}.
\]
The property-$j$ of this $\clp$ follows from the equivalence of \eqref{2.10} and \eqref{2.11}.
\end{Rk}
In view of Proposition \ref{PnW1}
we can define the Weyl function of the skew-self-adjoint Dirac system similar to the cases of the canonical system \cite{SaL3} and  self-adjoint Dirac system \cite{FKRSp1}. 
\begin{Dn} \label{defWeyl} The Weyl-Titchmarsh $($or simply Weyl$)$ function of Dirac system \eqref{1.1}, where   $\|v(x)\|$ is bounded  on $\BR_+$, is the function $\vp$
given by \eqref{2.3}.
\end{Dn}
\begin{Rk}\label{RkW1}
From Proposition  \ref{PnW1} we see that under condition \eqref{s1} the Weyl-Titchmarsh function always exists and that it is  unique as well.
\end{Rk}

Now we introduce a more general definition of the Weyl function.
For this reason we need the following corollary. 
\begin{Cy} \label{CyW2} Let the conditions of Proposition \ref{PnW1} hold.
Then the Weyl function is the unique function, which satisfies the inequality
\begin{align}&      \label{2.20}
\int_0^{\infty}
\begin{bmatrix}
I_{m_1} & \vp(z)^*
\end{bmatrix}
u(x,z)^*u(x,z)
\begin{bmatrix}
I_{m_1} \\ \vp(z)
\end{bmatrix}dx< \infty .
\end{align}
\end{Cy}
\begin{proof}.  Using the first inequality in \eqref{s2}  we derive
\begin{align}&      \label{2.21}
\int_0^{r}
\begin{bmatrix}
I_{m_1} & \vp(z)^*
\end{bmatrix}
u(x,z)^*u(x,z)
\begin{bmatrix}
I_{m_1} \\ \vp(z)
\end{bmatrix}dx \leq\frac{1}{2((\Im (z) - M)}\begin{bmatrix}
I_{m_1} & \vp(z)^*
\end{bmatrix}
\\ \nn &
\times \big(
\mf{A}(0,z)-\mf{A}(r,z)\big)
\begin{bmatrix}
I_{m_1} \\ \vp(z)
\end{bmatrix}.
\end{align}
Recall that \eqref{2.10} implies \eqref{2.11}, and so for $\vp$ satisfying \eqref{2.3} we get
\begin{align}\label{s4}
\begin{bmatrix}
I_{m_1} &  \vp(z)^*
\end{bmatrix}\mathfrak{A}(x,z)\begin{bmatrix}
I_{m_1} \\  \vp(z)
\end{bmatrix} \geq 0  \quad \mathrm{for} \,\,  \mathrm{all} \quad 0\leq x<\infty.
\end{align} 
Inequality \eqref{2.20} follows immediately from \eqref{2.21} and \eqref{s4}. 

Moreover,
since $u^*u\geq -\mf{A}$, inequality \eqref{2.12} yields
\begin{align}&      \label{2.22}
\int_0^{r}
\begin{bmatrix}0 &
I_{m_2} 
\end{bmatrix}
u(x,z)^*u(x,z)
\begin{bmatrix}
0 \\ I_{m_2} 
\end{bmatrix}dx \geq 
r I_{m_2} .
\end{align}
In view of \eqref{2.22}, the function $\vp$ satisfying \eqref{2.20} is unique.
\end{proof}
\begin{Dn} \label{defWeyl2} The Weyl function of Dirac system \eqref{1.1}  on $[0, \, \infty)$,
where $\|v\|$ is locally summable, is the function $\vp$
satisfying \eqref{2.20}.
\end{Dn}
From  Corollary \ref{CyW2}, we see that Definition \ref{defWeyl2}
is equivalent to Definition \ref{defWeyl} for the case that the potentials
are bounded, but it can also be used in more general situations \cite{CG2}.
A definition of the form \eqref{2.20} is a more classical one and deals with
solutions of \eqref{1.1} which belong to $L^2(0, \, \infty)$. Compare 
\eqref{2.20} with definitions
of the Weyl-Titchmarsh or $M$-functions for discrete and continuous
systems in \cite{LS, Mar, SaL3, T1, T2} (see also references therein).
For the case of system \eqref{1.1}, where $m_1=m_2$, 
Definition \ref{defWeyl2} was introduced and the corresponding Weyl function was studied
in \cite{CG2, GKS2, SaA1, SaA2, ALS06}.
\section{The inverse problem on a finite interval} \label{fi}
\setcounter{equation}{0}
In this section we consider a procedure for solving the inverse problem,
that is, for recovering the potential $v$ (and, hence, the Dirac system) from the
Weyl function. We note that it suffices to recover $v$ almost everywhere on all the finite intervals $[0, \,l]$. Therefore, we introduce the corresponding definition of  Weyl functions
of the skew-self-adjoint Dirac system on an arbitrary fixed interval $[0, \, l]$ and solve the inverse
problem on this interval.
\begin{Dn} \label{defWeyll} Weyl-Titchmarsh $($or simply Weyl$)$ functions of Dirac system \eqref{1.1} on $[0, \, l]$ which satisfies \eqref{s1'},  are the  functions of the form
\eqref{2.2}, where $x=l$
and $\clp$
satisfies \eqref{2.2'}, i.e., the analytic functions $\vp(z)\in \cln(l,z)$ $\, (z \in \BC_M)$.
\end{Dn}
Recall that in view of Remark \ref{Rk3} the set $ \cln(l,z)$ is well-defined.

We assume that $v$  satisfies \eqref{s1'} and let
\begin{align}&      \label{3.1}
\b(x)=\begin{bmatrix}
I_{m_1} & 0
\end{bmatrix}u(x,0), \quad \g(x)=\begin{bmatrix}
0 &I_{m_2}
\end{bmatrix}u(x,0).
\end{align}
It follows from $\sup_{x < l}\|v\|<\infty$ and from \eqref{1.1} that
\begin{align}&      \label{3.2}
\sup_{x < l}\|\g^{\prime}(x)\|<\infty, \quad \g^{\prime}:=\frac{d}{dx}\g .
\end{align}
Moreover, from \eqref{1.1}-\eqref{1.3} we get 
\begin{align}&      \label{3.2'}
u(x,\ov z)^*u(x,z)=u(x,z)u(x, \ov z)^*=I_m.
\end{align}
Therefore, 
\eqref{3.1} implies
\begin{align}&      \label{3.3}
\b \b^*\equiv I_{m_1}, \quad \g  \g^*\equiv I_{m_2}, \quad \b \g^*\equiv 0. 
\end{align}
Next, as in the self-adjoint case \cite{FKRSp2}, we need the  similarity result below for the Volterra operator
\begin{align}&      \label{3.4}
K=\int_0^x  \, F(x)G(t) \cdot \, dt, \quad K \in B\big(L^2_{m_1}(0, \, l)\big),
\end{align}
where $F(x)$ is an $m_1 \times m$ matrix function, $G(t)$ is an $m \times m_1$ matrix function,
and
\begin{equation} \label{3.5}
\I F(x)G(x) \equiv I_{m_2}.
\end{equation}
\begin{Pn}\label{PnSim}{\rm \cite{SaL0}} Let $F$ and $G$ be boundedly differentiable
and let \eqref{3.5} hold. Then we have
\begin{equation} \label{3.6}
K= EAE^{-1}, \quad A:=-\I \int_0^x \, \cdot \, dt,  \quad A, E, E^{-1} \in B\big(L^2_{m_2}(0, \, l)\big),
\end{equation}
where $K$ is given by \eqref{3.4} and $E$ is a triangular operator of the form
\begin{equation} \label{3.7}
(E f)(x)=\rho(x) f(x)+\int_0^x E(x,t)f(t)dt, \quad \frac{d}{d x}\rho =\I F^{\prime}G\rho, \quad
\det \rho(0)\not=0.
\end{equation}
Moreover, the operators $E^{\pm 1}$ map functions with bounded derivatives
into functions with bounded derivatives.
\end{Pn}
We set
\begin{align}&      \label{3.8}
F=\g, \quad G=-\I \g^*, \quad \g \g^*\equiv I_{m_2}, \quad \g \in B^1[0,\,l],
\end{align}
where $B^1$ stands for the class of boundedly differentiable matrix functions. 
Clearly $F$ and $G$ in \eqref{3.8} satisfy the conditions of Proposition \ref{PnSim}. 

We partition $\g(x)$ into two blocks 
$\g=\begin{bmatrix}
\g_1 & \g_2
\end{bmatrix}$, where $\g_1, \, \g_2$ are,
respectively, $m_2\times m_1$ and $m_2\times m_2$ matrix functions. 
We show that, without loss of generality, one can choose $E$ satisfying condition
\begin{align}&      \label{3.12}
E^{-1}\g_2\equiv I_{m_2},
\end{align}
where $E^{-1}$ is applied to $\g_2$ columnwise.
\begin{Pn} \label{PnSimN} Let 
$K$ be given by \eqref{3.4}, where $F$ and $G$ satisfy \eqref{3.8}, and let $\wt E$
be a similarity operator from Proposition \ref{PnSim}.
Suppose $E_0 \in B\big(L^2_{m_2}(0, \, l)\big)$ is determined by the equalities
\begin{align}&      \label{3.13}
\big(E_0f\big)(x)=\rho(0)^{-1}\g_2(0)f (x) + \int_0^x E_0(x-t) f(t) dt, \quad 
E_0(x):=\big(\wt E^{-1}\g_2\big)^{\prime}(x).
\end{align}
Then the operator $E:=\wt E E_0$ is another similarity operator from Proposition \ref{PnSim},
which satisfies the additional condition   \eqref{3.12}.
\end{Pn}
\begin{proof}.
The proof of the proposition is similar to the case $m_1=m_2$ (see, e.g.,
\cite[pp. 103, 104]{SaL0}) and coincides with the proof of Proposition 2.2 in \cite{FKRSp2}.
\end{proof}
\begin{Rk}\label{kern} The kernels of the operators $\wt E^{\pm 1}$, which are
constructed in \cite{SaL0}, as well as the kernels of the operators $E_0^{\pm 1}$
from Proposition \ref{PnSimN} are bounded.
Therefore, without loss of generality, we always assume that the kernels of $E^{\pm 1}$ are bounded.
\end{Rk}
The next lemma easily follows from Proposition \ref{PnSimN}
and will be used to construct a representation of the fundamental solution $u$.
\begin{La}\label{LaNode} Let $\g$ be an $m_2 \times m$ matrix function, which satisfies the last two relations in \eqref{3.8}, and set
\begin{align}&      \label{p6}
S:=E^{-1}\big(E^*\big)^{-1}, \quad \Pi:= \begin{bmatrix}
\Phi_1 & \Phi_2
\end{bmatrix}, \quad 
\Phi_k \in B\big(\BC^{m_k}, \, L^2_{m_2}(0, \, l)\big);
\\
&      \label{3.19}
\big(\Phi_1 g\big)(x)=\Phi_1(x)g, \quad
\Phi_1(x):=\big(E^{-1}\g_1\big)(x);  \quad  \Phi_2 g=I_{m_2}g\equiv g;
\end{align}
where $E$ is constructed $($for the given $\g)$ in Proposition \ref{PnSimN} and
$g$ is used to denote finite-dimensional vectors.
Then $A$, $S$, and $\Pi$ form an $S$-node, that is $($see \cite{SaL1, SaL2, SaL3}$)$,
the operator identity
\begin{align}&      \label{p9}
AS-SA^*=-\I \Pi  \Pi^*
\end{align}
holds. Furthermore,
we have
\begin{align}&      \label{p7}
\ov{\sum_{i=0}^{\infty}\im \Big(\big(A^*\big)^i S^{-1}\Pi\Big)}=L^2_{m_2}(0, \, l).
\end{align}
\end{La}
\begin{proof}.  Because of \eqref{3.4},  \eqref{3.6} and \eqref{3.8}, we get
\begin{align}&      \label{p8}
EAE^{-1}-\big(E^{-1}\big)^*A^*E^*=K-K^*=-\I\g(x)\int_0^l\g(t)^*\,\cdot \,dt.
\end{align}
The second equality in \eqref{p6}
and formulas \eqref{3.12} and \eqref{3.19} lead us to the equalities
\begin{align}&      \label{p10}
\Pi g=\big(E^{-1}\g\big)(x)g, \quad \Pi^*f=\int_0^l \big(E^{-1}\g\big)(t)^*f(t)dt=
\int_0^l \g(t)^*\big((E^{-1})^*f\big)(t)dt.
\end{align}
Using the first equality in \eqref{p6} and relations \eqref{p10} we can rewrite \eqref{p8}
as the operator identity \eqref{p9}.
To prove \eqref{p7} we will show that
\begin{align}&      \label{p11}
\sum_{i=0}^{N}\im \Big(\big(A^*\big)^i S^{-1}\Pi\Big)
\supseteq \sum_{i=0}^{N}\im \big(S^{-1}A^i\Pi\big)
=S^{-1}\sum_{i=0}^{N}\im \big(A^i\Pi\big).
\end{align}
For this reason, we rewrite \eqref{p9} as $S^{-1}A=A^*S^{-1}-\I
S^{-1}\Pi \Pi^* S^{-1}$. Hence, for $N_1, N_2 \geq 0$ we obtain
\begin{align}&     \label{p12}
\im \Big(\big(A^*\big)^{N_1+1} S^{-1}A^{N_2}\Pi\Big)
\\ \nn &
+
\sum_{i=0}^{N_1+N_2}\im \Big(\big(A^*\big)^i S^{-1}\Pi\Big)
\supseteq \im \Big(\big(A^*\big)^{N_1} S^{-1}A^{N_2+1}\Pi\Big).
\end{align}
Using \eqref{p12}, we derive \eqref{p11} by induction.
In view of \eqref{p11}, it suffices to show that
\begin{align}&      \label{p13}
\ov{\sum_{i=0}^{\infty}\im \big(A^i \Pi\big)}=L^2_{m_2}(0, \, l),
\end{align}
which, in turn, follows from 
$ \Pi= \begin{bmatrix}
\Phi_1 & \Phi_2
\end{bmatrix}$, where
$\Phi_2 g \equiv g$ (see \eqref{3.19}).
\end{proof}
\begin{Rk} \label{TrMF}
Given an $S$-node \eqref{p9}, we introduce a transfer matrix function
in Lev Sakhnovich form $($see \cite{SaL1, SaL2, SaL3}$):$
\begin{align}&      \label{p13a}
w_A(r,z):=I_m-\I z\Pi^*S_r^{-1}(I-zA_r)^{-1}P_r\Pi, \quad 0<r\leq l,
\end{align}
where $I$ is the identity operator; $A_r, \, S_r \, \in \, B\big(L^2_{m_2}(0, \, r)\big)$, 
\begin{align}&      \label{p13b}
A_r:=P_rAP_r^*, \quad S_r:=P_rSP_r^*;
\end{align}
 $A$ is given by \eqref{3.6}, the operators $S$ and $\Pi$ are given by
 \eqref{p6} and \eqref{3.19}, and the
operator $P_r$ is an orthoprojector from $L^2_{m_2}(0, \, l)$ on $L^2_{m_2}(0, \, r)$ such that
\begin{align}&      \label{p13c}
\big(P_rf\big)(x)=f(x) \quad (0<x<r), \quad f \in L^2_{m_2}(0, \, l).
\end{align}
Since $P_rA=P_rAP_r^*P_r$, it follows from \eqref{p9} that the operators $A_r$, $S_r$, and $P_r\Pi$ also form an 
$S$-node, i.e., the operator identities
\begin{align}&      \label{p9'}
A_rS_r-S_rA_r^*=-\I P_r \Pi  \Pi^*P_r^*
\end{align}
hold.
\end{Rk}
Now, in a way similar to \cite{SaA0, SaA1} the fundamental solution $w$
of the system
\begin{align}&      \label{3.14}
\frac{d}{dx}w(x,z)=-\I z \g(x)^*\g(x)w(x,z), \quad w(0,z)=I_m
\end{align}
is constructed.
\begin{Tm}\label{FundSol}
Let $\g$ be an $m_2 \times m$ matrix function, which satisfies the last two relations in \eqref{3.8}.
 Then, the fundamental solution $w$ given by \eqref{3.14} 
 admits the representation
 \begin{align}&      \label{3.16}
w(r,z)=w_A(r,z),
\end{align}
where $w_A(r,z)$ is defined in Remark \ref{TrMF} $($see \eqref{p13a}$)$.
\end{Tm}
\begin{proof}. The statement of the theorem follows from the Continual Factorization Theorem
(see \cite[p. 40]{SaL3}). More precisely, our statement follows from a corollary of the
Continual Factorization Theorem, namely, from Theorem 1.2 \cite[p. 42]{SaL3}.
Using Lemma \ref{LaNode} we easily check that the conditions 
of Theorem 1.2 \cite[p. 42]{SaL3} are fulfilled. (We note that, though \eqref{p7} is valid,
 this condition  is not essential for the factorization  Theorem 1.2 \cite[p. 42]{SaL3} because $S$ is invertible.)
Theorem 1.2 \cite{SaL3} implies: if $\Pi^*S_r^{-1}P_r\Pi$
is boundedly differentiable, then we have
\begin{align}&      \label{p14}
\frac{d}{dr}w_A(r,z)=-\I zH(r)w_A(r,z), \quad \lim_{r\to +0}w_A(r,z)=I_m, \\
&      \label{p14'}
 H(r):=\frac{d}{dr}\big(\Pi^*S_r^{-1}P_r\Pi\big),
\end{align} 
where $w_A$ is given by \eqref{p13a}. 
Taking into account that $E^{\pm 1}$ are lower triangular operators, we derive
\begin{align}&      \label{p14!}
P_rEP_r^*P_r=P_rE, \qquad  \big(E^{-1}\big)^*P_r^*=P_r^*P_r\big(E^{-1}\big)^*P_r^*.
\end{align} 
Hence, formulas \eqref{p6} and
\eqref{p13b} lead us to
\begin{align}&      \label{p14!!}
 S_r^{-1}=E_r^*E_r, \qquad E_r:=P_rEP_r^*.
\end{align} 
Therefore, using  \eqref{p10}, we rewrite  \eqref{p14'} as
\begin{align}&      \label{p15}
 H(r)=\g(r)^*\g(r).
\end{align} 
Compare formulas \eqref{3.14} and
\eqref{p14}, \eqref{p15}  to see that \eqref{3.16} holds.
\end{proof}
Now, consider again the case of a Dirac system.  Because of  \eqref{3.1} and \eqref{3.3} we obtain
\begin{align}&      \label{p16}
u(x,0)\g(x)^*\g(x)u(x,0)^{-1}=\begin{bmatrix}
0 & 0\\0 & I_{m_2}
\end{bmatrix}.
\end{align}
In view of \eqref{3.14} and \eqref{p16}, direct differentiation shows that the following corollary of Theorem \ref{FundSol} holds.
\begin{Cy}\label{FSD} Let $u(x,z)$ be the fundamental solution of a 
skew-self-adjoint Dirac system, where the potential
$v$ is bounded, and let $\g$ be given by \eqref{3.1}. Then $u(x,z)$ admits the representation
\begin{align}&      \label{p17}
u(x,z)=\E^{ixz}u(x,0)w(x,2z),
\end{align}
 where $w$ has the form \eqref{3.16} and the $S$-node generating the transfer matrix function $w_A$    
is recovered from $\g$ in Lemma \ref{LaNode}.
\end{Cy}
\begin{Rk}\label{beta'}
For the case that $\g$ is given by \eqref{3.1}, it follows from   \eqref{1.1} and \eqref{3.3} that
\begin{align}&      \label{3.9}
\g^{\prime}(x)\g(x)^*=- \begin{bmatrix}
v(x)^* & 0
\end{bmatrix}u(x,0)\g(x)^*=- v(x)^*\b(x)\g(x)^*\equiv 0.
\end{align}
Thus, from \eqref{3.8} and \eqref{3.9} we see that 
$F^{\prime}G\equiv 0$, and so $\rho$ in \eqref{3.7} is a constant matrix.
Therefore, since $\g_2(0)=I_{m_2}$ and $\Phi_2(x) \equiv I_{m_2}$, equality \eqref{p10} implies that  $\rho(x)\equiv I_{m_2}$, and 
formula \eqref{3.7} can be rewritten in the form
\begin{equation} \label{3.10}
(E f)(x)=f(x)+\int_0^x E(x,t)f(t)dt, \quad E \in B\big(L^2_{m_2}(0, \, l)\big).
\end{equation}
\end{Rk}

Important publications by  F. Gesztesy and B. Simon \cite{GeSi0, GeSi, Si}  
gave rise to a whole series of interesting papers and results
on the high energy asymptotics of Weyl functions and  local Borg-Marchenko-type
uniqueness theorems (see, e.g., \cite{CG1, CGZ, CGZ2, FKRSp2, FKS2, LaWo, 
SaA3, ALS06} and references therein). Here we use the $S$-node scheme
from \cite{SaA1} and generalize the high energy
asymptotics result from \cite{ALS06} for the case that Dirac system \eqref{1.1}
has  a rectangular $m_1 \times m_2$ potential $v$, where $m_1$ is not necessarily equal to $m_2$. 
We  first recall that $S>0$ and $\Phi_1(x)=\big(E^{-1}\g_1\big)(x)$ is boundedly differentiable in Lemma \ref{LaNode}. Next, we note that since $\g_1(0)=0$ and
$E$ has the form \eqref{3.10}, the equality
\begin{align} \label{S1}
\Phi_1(0)=0
\end{align}
holds.
Therefore, using   Theorem 3.1 from \cite{FKRSp2}  and substituting \eqref{S1} there,
we get the statement below.
\begin{Tm}\label{TmIdent}
Let  $\Pi=[\Phi_1 \quad \Phi_2]$ be constructed in Lemma \ref{LaNode}.
Then there is a unique solution 
$S \in B\big(L^2_{m_2}(0,l)\big)$ of the operator
identity \eqref{p9}. This $S$ is strictly positive $($i.e., $S>0)$ and is  defined by the equalities
\begin{align}&
\label{4.5}
S= I +\int_0^ls(x,t)\, \cdot \, dt, \quad s(x,t):=\int_0^{\min(x,t)}\Phi_1^{\prime}(x-\zeta)\Phi_1^{\prime}(t-\zeta)^*d\zeta.
\end{align}
\end{Tm}
Now, we can apply the $S$-node scheme  
to derive the high energy asymptotics. 
\begin{Tm}\label{TmHea} Assume that $\vp $ is a Weyl function $($i.e., $\vp$ admits representation \eqref{2.2}, where $x=l$ and $\clp$
satisfies \eqref{2.2'}$)$,
 and the potential $v$ of the corresponding   Dirac system \eqref{1.1} 
is bounded on $[0, \, l]$ $($i.e., \eqref{s1'} holds$)$. Then $($uniformly with respect to $\Re(z))$ we have
\begin{align} \label{hea}&
\vp(z)=2\I z\int_0^l\E^{2\I xz}\Phi_1(x)dx+O\left(2z \E^{2\I lz}/ \sqrt{\Im(z)}\right), \quad \Im (z) \to \infty.
\end{align}
\end{Tm}
\begin{proof}. {\bf Step 1.} Similar to the self-adjoint case,  we first show that
\begin{align}&\label{4.11}
\left\| (I-2zA)^{-1}\Pi 
\begin{bmatrix}
I_{m_1} \\ \vp(z)
\end{bmatrix}\right\| \leq C/\sqrt{\Im z} \quad {\mathrm{for}}\,\, 
{\mathrm{some}} \quad C>0.
\end{align}
For this purpose, we consider the matrix function
\begin{align}&\label{4.6}
\clu(z)=\begin{bmatrix}
I_{m_1} & \vp(z)^*
\end{bmatrix}\big(
w_A(l,2z)^*w_A(l,2z)-I_m\big)
\begin{bmatrix}
I_{m_1} \\ \vp(z)
\end{bmatrix}.
\end{align}
The proof that this function is bounded is different from the self-adjoint case (compare
with \cite{FKRSp2}). More specifically,
because of \eqref{3.2'}, \eqref{3.16} and \eqref{p17} we have
\begin{align}&\nn
\clu(z)=e^{\I l(\ov z-  z)}\begin{bmatrix}
I_{m_1} & \vp(z)^*
\end{bmatrix}
u(l,z)^*u(l,z)
\begin{bmatrix}
I_{m_1} \\ \vp(z)
\end{bmatrix}
-I_{m_1}-\vp(z)^*\vp(z).
\end{align}
We  substitute \eqref{2.9} into the formula above and rewrite it as
\begin{align}\nn 
\clu(z)=&e^{\I l(\ov z-  z)}
\Big(\big(\begin{bmatrix}
I_{m_1} & 0
\end{bmatrix}u(l,z)^{-1}\clp(z)\clp_1(z)^{-1}\big)^{-1}\Big)^*
\big(\clp(z)\clp_1(z)^{-1}\big)^*
\\& \nn \times
\big(\clp(z)\clp_1(z)^{-1}\big)
\big(\begin{bmatrix}
I_{m_1} & 0
\end{bmatrix}u(l,z)^{-1}\clp(z)\clp_1(z)^{-1}\big)^{-1}
\\& \label{4.8}
 -I_{m_1}-\vp(z)^*\vp(z),
\end{align}
where $\clp_1$ is the upper $m_1\times m_1$ block of $\clp$. Clearly,
inequalities \eqref{2.1} imply the invertibility of  $\clp_1(z)$ and boundedness of $\clp(z)\clp_1(z)^{-1}$. 

{\bf Step 2.} To derive \eqref{4.11} we also need to examine the asymptotics of $u(l,z)$ ($z\to \infty$).
This is achieved by considering the procedure for constructing the transformation
operator for Dirac system (see also some more complicated cases in \cite{KoSaTe0, MST, SaL0}).
We will show that $u$ admits an integral representation
\begin{align} \label{intre}&
u(x,z)=\E^{\I x z j}+\int_{-x}^{x}\E^{\I t z}N(x,t)dt, \quad \sup \|N(x,t)\|<\infty \quad
(0<|t|<x<l).
\end{align}
Indeed, it is easily checked that
\begin{align} \label{IR1}&
u(x,z)=\sum_{i=0}^{\infty}\nu_i(x,z), 
\\ & \nn
 \nu_0(x,z):=\E^{\I x z j}, \quad \nu_i(x,z):=\int_0^x \E^{\I (x-t) z j}jV(t)\nu_{i-1}(t,z)dt \quad \mathrm{for} \quad i>0,
\end{align}
since the right-hand side of \eqref{IR1} satisfies both equation \eqref{1.1} and the normalization condition at $x=0$. Moreover, using induction we derive in a standard way that
\begin{align} \label{IR2}&
\nu_i(x,z)=\int_{-x}^{x}\E^{\I t z}N_{i}(x,t)dt, \quad \sup \|N_i(x,t)\| \leq (2M)^ix^{i-1}/(i-1)!,
\end{align}
where $i>0$ and $M$ is the same as in \eqref{s1'}. We see that formulas \eqref{IR1} and 
\eqref{IR2}
yield \eqref{intre}.

{\bf Step 3.}  Because of \eqref{3.2'} and \eqref{intre} we get
\begin{align} \label{IR3}&
u(l,z)^{-1}=u(l, \ov z)^*=\E^{-\I l z}\left( \begin{bmatrix}
I_{m_1} & 0 \\
0 & \E^{2 \I l z}I_{m_2}
\end{bmatrix} +o(1) \right) \quad (z\in \BC_+, \,\, z \to \infty)
.
\end{align}
In view of \eqref{2.2}, \eqref{4.8} and \eqref{IR3}, for any $\ve >0$ there are  numbers
$\wt M>M>0$ and $\wt C>0$, such that for all $z$ satisfying $\Im z>\wt M$ we have
\begin{align} \label{IR4}&
\|\vp(z)\| \leq \ve, \quad \|\clu(z)\| \leq \wt C.
\end{align}

Next, using  \eqref{p9} and  \eqref{p13a}
we easily obtain  (see, e.g., \cite{SaL2, FKS2})
\begin{align} \label{4.7}&
w_A(l, z)^*w_A(l,z)=I_m+\I(\ov z -z)\Pi^*(I-\ov z A^*)^{-1}S_l^{-1}(I-zA)^{-1}\Pi .
\end{align}
We substitute \eqref{4.7} into  \eqref{4.6} to rewrite the second inequality in \eqref{IR4} 
in the form 
\begin{align}&\label{4.10}
 2\I (\ov z-z) \begin{bmatrix}
I_{m_1} & \vp(z)^*
\end{bmatrix}
\Pi^*(I-2\ov z A^*)^{-1}S_l^{-1}(I-2zA)^{-1}\Pi 
\begin{bmatrix}
I_{m_1} \\ \vp(z)
\end{bmatrix}\leq \wt C I_{m_1}.
\end{align}
Since $S$ is invertible, positive and bounded, inequality \eqref{4.10} yields \eqref{4.11}.

{\bf Step 4.}
We easily check directly (see also these formulas in the works on the case
$m_1=m_2$) that
\begin{align}  \label{4.12.0}&
(I-zA)^{-1}f=f-\I z \int_0^x\E^{\I(t-x)z}f(t)
dt,
\\ \label{4.12}&
\Phi_2^*(I-2zA)^{-1}f=\int_0^l\E^{2\I(x-l)z}f(x)dx, 
\\  \label{4.13}&
  \Phi_2^*(I-2zA)^{-1}\Phi_2=\frac{\I}{2z}\big(\E^{-2\I lz}-1\big)I_{m_2}.
\end{align}
Now, apply $-\I\Phi_2^*$ to the operator on the
left-hand side of \eqref{4.11}  and use \eqref{4.12} and  \eqref{4.13} to get
(uniformly with respect to $\Re(z)$) the equality
\begin{align} \label{4.14}&
\frac{1}{2z}\big(\E^{-2\I lz}-1\big)\vp(z)=\I \E^{-2\I lz}\int_0^l\E^{2\I xz}\Phi_1(x)dx+
O\left(\frac{1}{\sqrt{\Im(z)}}\right).
\end{align}
Because of  \eqref{4.14} (and the first inequality in \eqref{IR4}), we see   that \eqref{hea} holds.
\end{proof}
 The  integral representation below follows from the high-energy asymptotics of $\vp$ and
is essential in interpolation and inverse problems.
\begin{Cy}\label{cyHea}  Let $\vp$ be a Weyl function of a Dirac system on $[0, \, l]$ which
satisfies the condition \eqref{s1'} of Theorem \ref{TmHea}.
Then we have
\begin{align} \label{5.1}&
\Phi_1\Big(\frac{x}{2}\Big)=\frac{1}{\pi}\E^{x\eta}{\mathrm{l.i.m.}}_{a \to \infty}
\int_{-a}^a\E^{-\I x\xi}\frac{\vp(\xi+\I \eta)}{2\I(\xi +\I \eta)}d\xi, \quad \eta >M,
\end{align}
where l.i.m. stands for the entrywise limit in the norm of  $L^2(0,2l)$.
\end{Cy}
\begin{proof}. According to Definition \ref{defWeyll}  and relations \eqref{2.16} and \eqref{2.17},
the matrix function $\vp$ is non-expansive in $\BC_M$.
Since $\vp$ is analytic and non-expansive in $\BC_M$,
 it admits (see, e.g., \cite[Theorem V]{WP}) a representation
\begin{align} \label{repr0}&
\vp(z)=2\I z\int_0^{\infty}\E^{2\I xz}\Phi (x)dx, \quad z=\xi+\I \eta, \quad \eta>M>0,
\end{align}
where $\E^{-2\eta x}\Phi(x) \in L^2_{m_2\times m_1}(0, \, \infty)$.
Because of \eqref{hea} and \eqref{repr0} we derive
\begin{align}&\nn
\int_0^l\E^{2\I (x-l)z}(\Phi_1(x)-\Phi(x))dx=\int_l^{\infty}\E^{2\I (x-l)z}\Phi(x)dx+O(1/\sqrt{\Im(z)})
\end{align}
for $ \Im(z) \to +\infty $. Taking into account that $\E^{-2(M+\ve) x}\Phi(x) \in L^2_{m_2\times m_1}(0, \, \infty)$ for $\ve >0$, we rewrite the formula above as
\begin{align}&\label{IR5}
\int_0^l\E^{2\I (x-l)z}(\Phi_1(x)-\Phi(x))dx=O(1/\sqrt{\Im(z)}),
\end{align}
and the equality \eqref{IR5} is uniform with respect to $\Re(z)$.
Clearly, the left-hand side of \eqref{IR5} is bounded in the domains, where $\Im (z)$
is bounded from above. Hence, in view of  \eqref{IR5}, its left-hand side is bounded
in all of $\BC$ and tends to zero on some rays. The identities
\begin{align}&\label{IR6}
\int_0^l\E^{2\I (x-l)z}(\Phi_1(x)-\Phi(x))dx \equiv 0, \quad \mathrm{i.e.,} 
\quad \Phi_1(x)\equiv \Phi(x)
\end{align}
are clear. 

Using the Plancherel Theorem,
we apply the inverse Fourier transform to formula \eqref{repr0} and get a representation
of the form \eqref{5.1} for $\Phi$. Since, according to \eqref{IR6}, we have $ \Phi_1\equiv \Phi$,
the formula  \eqref{5.1} is valid.
\end{proof}
Directly from \eqref{1.1} and \eqref{3.1}, we easily get a useful formula
\begin{align} \label{5.4'}&
\b^{\prime}(x)= v(x)\g(x),
\end{align}
which, because of \eqref{3.3},  implies
\begin{align} \label{5.5}&
\b^{\prime}\b^*=0, \qquad \b^{\prime}\g^*= v.
\end{align}
Thus, we will solve the inverse problem and recover $v$ from $\vp$ once we recover
$\b$ and $\g$. The recovery of $\Phi_1$ from $\vp$ is studied in Corollary \ref{cyHea}.
The next step to solve the inverse problem is to recover $\b$ from $\Phi_1$ in the
following proposition.
\begin{Pn} \label{Pnbeta} Let Dirac system \eqref{1.1} on $[0, \, l]$ satisfy
the conditions of Theorem \ref{TmHea} and let $\Pi$ and $S$ be the operators
constructed in Lemma \ref{LaNode}. Then the matrix function $\b$, which is defined in \eqref{3.1},
satisfies the equality 
\begin{align}&      \label{5.9}
\b(x)=\begin{bmatrix}I_{m_1} &0 \end{bmatrix}
-\int_0^x\Big(S_x^{-1}\Phi_1^{\prime}\Big)(t)^*\begin{bmatrix}\Phi_1(t) & I_{m_2} \end{bmatrix}dt.
\end{align} 
 \end{Pn}
\begin{proof}. First, we use \eqref{3.12} and \eqref{3.19} to  express $\g$ in the form 
 \begin{align}&      \label{5.10}
\g(x)=\big(E \begin{bmatrix}\Phi_1 & I_{m_2} \end{bmatrix}\big)(x).
\end{align} 
It follows that
 \begin{align}&      \label{5.11}
\I \big(E AE^{-1}E\Phi_1^{\prime}\big)(x)=\g_1(x)-\g_2(x)\Phi_1(+0).
\end{align}
According to \eqref{3.1}, \eqref{3.10}, \eqref{5.10} and Remark \ref{kern},
we have 
 \begin{align}&      \label{IR7}
\g_1(0)=0, \quad \Phi_1(+0)=0.
\end{align}
Therefore, we rewrite \eqref{5.11} as
 \begin{align}&      \label{5.12}
\I \big(E AE^{-1}E\Phi_1^{\prime}\big)(x)=\g_1(x).
\end{align}
Next, we substitute $K=EAE^{-1}$ from \eqref{3.6} into \eqref{5.12} and (using  \eqref{3.4} and \eqref{3.8}) we  get
 \begin{align}&      \label{5.13}
\g_1(x)=\g(x)\int_0^x \g(t)^*\Big(E\Phi_1^{\prime}\Big)(t)dt.
\end{align}
Formulas \eqref{5.10} and \eqref{5.13} imply
 \begin{align}&      \label{5.14}
\g_1(x)=\g(x)\int_0^x \big(E \begin{bmatrix}\Phi_1 & I_{m_2} \end{bmatrix}\big)(t)^*\Big(E\Phi_1^{\prime}\Big)(t)dt.
\end{align}
Because of  \eqref{p14!}, \eqref{p14!!} and \eqref{5.14} we see that
\begin{align} &    \label{5.15}
\g(x)\t(x)^*\equiv 0 \quad \mathrm{for} \quad \t(x):=\begin{bmatrix}I_{m_1} &0 \end{bmatrix}
-\int_0^x 
\Big(E\Phi_1^{\prime}\Big)(t)^*
\\ & \nn \times
\big(E \begin{bmatrix}\Phi_1 & I_{m_2} \end{bmatrix}\big)(t)dt 
 =
\begin{bmatrix}I_{m_1} &0 \end{bmatrix}
-\int_0^x\Big(S_x^{-1}\Phi_1^{\prime}\Big)(t)^*\begin{bmatrix}\Phi_1(t) & I_{m_2} 
\end{bmatrix}dt,
\end{align}
where $S_x:=P_xSP_x^*$. We will show that $\t=\b$.

In view of \eqref{5.10} and the second relation in \eqref{5.15},
we have 
\begin{align}     \label{5.15'}&
\t^{\prime}(x)=-\Big(E\Phi_1^{\prime}\Big)(x)^*\g(x).
\end{align}
Therefore, \eqref{3.3}
leads us to
\begin{align}     \label{5.16}&
\b(x)\t^{\prime}(x)^* \equiv 0.
\end{align}
We furthermore compare \eqref{3.3} with the first  relation in \eqref{5.15}
to see that
\begin{align}     \label{5.17}&
\t(x)=\vk(x)\b(x),
\end{align}
where $\vk(x)$ is an $m_1 \times m_1$ matrix function.
According to \eqref{3.1} and \eqref{5.15'}, respectively, $\b$ and $\t$ are boundedly
differentiable on $[0, \, l]$, and so $\vk$
  is also boundedly differentiable.
Now, equalities \eqref{5.16}, \eqref{5.17}
and the first relations in \eqref{3.3}, \eqref{5.5} yield that $\vk^{\prime}\equiv 0$
(i.e., $\vk$ is a constant). It follows from \eqref{3.1}, \eqref{5.15} and \eqref{5.17}
that $\vk(0)=I_{m_1}$, and therefore $\vk \equiv I_{m_1}$, that is, $\t \equiv \b$.
Thus, \eqref{5.9} follows directly from \eqref{5.15}.
\end{proof}
In a way which is quite similar to the proof of  \eqref{5.4'} and \eqref{5.5} we get
\begin{align}     \label{IR8}&
\g^{\prime}=-v^*\b, \quad \g^{\prime}\g^*=0.
\end{align}
In order to recover $\g$ from $\b$, we recall that $\b$ is boundedly differentiable and $\b\b^* \equiv I_{m_1}$. Therefore, an $m_2 \times m$
matrix function $\wt \g$, which is continuous, peacewise  differentiable (with the differentiability dissappearing, possibly, on a finite
set only), having bounded left and right derivatives everywhere on $[0,\, l]$  and satisfying relations
\begin{align}     \label{IR9}&
\b \wt \g^*\equiv 0, \quad \wt \g \wt \g^*>0, \quad \wt \g(0)=\g(0),
\end{align}
can be easily constructed. 
Since
\begin{align}     \label{IR10}&
\b\g^*=\b \wt \g^*= 0, \quad \g\g^*=I_{m_2},  \quad \wt \g \wt \g^*>0, \quad \wt \g(0)=\g(0),
\end{align}
the matrix function $\g$ admits representation
\begin{align}     \label{IR11}&
\g=\wt \vk \wt \g, \quad \wt \vk \wt \vk^*>0, \quad \wt\vk(0)=I_{m_2},
\end{align}
where $\wt \vk$ is continuous and peacewise differentiable. In view of \eqref{IR11}, the second equality in \eqref{IR8} can be rewritten in the form
\begin{align}     \label{IR12}&
\wt \vk^{\prime}=-\wt \vk\wt \g^{\prime} \wt \g^* (\wt \g \wt \g^*)^{-1}, \quad \wt\vk(0)=I_{m_2},
\end{align}
which uniquely defines $\wt \vk$. Hence, taking into account Corollary \ref{cyHea}, Proposition 
\ref{Pnbeta} and formula \eqref{5.5}, we see that we have obtained a procedure to solve the inverse problem.
\begin{Tm}\label{TmIP} Let the potential $v$ of the  Dirac system \eqref{1.1} 
 on $[0, \, l]$ be bounded $($i.e., let  \eqref{s1'} hold$)$. Then $v$ can be uniquely
 recovered from the Weyl function $\vp$ in the following way.
 First, we recover $\Phi_1$ using formula \eqref{5.1}. Next, we use formula \eqref{5.9}, where 
 the operator $S_x=P_xSP_x^*$ and $S$
 is given by \eqref{4.5}, to recover $\b$. We recover $($from $\b)$ matrix functions $\wt \g$
 and $\wt \vk$ via formulas \eqref{IR9} and \eqref{IR12}, respectively. The matrix function
 $\g$ satisfies the equality $\g=\wt \vk \wt \g$. Finally, we get $v$ from the equality
\begin{align}     \label{IR13}&
v(x)=\b^{\prime}(x)\g(x)^*.
\end{align}
\end{Tm}
The last statement in this section is a Borg-Marchenko-type uniqueness
theorem, which follows from Theorems \ref{TmHea} and  \ref{TmIP}.
More precisely, we will discuss the Weyl functions $\vp$ and $\wh \vp$ of two Dirac systems.
Matrix functions associated with $\wh \vp$ are written with a ``hat'' (e.g.,
$\wh v, \,\wh \Phi_1$).
\begin{Tm}\label{BM} Let $\vp$ and $\wh \vp$ be  Weyl functions
of two Dirac systems on $[0, \, l]$ with  bounded potentials,
denoted, respectively, by $v$ and $\wh v$. In other words, we assume that 
\begin{align}     \label{IR14}&
\max (\|v(x)\|,
\, \|\wh v(x)\|) \leq M, \quad
0<x<l.
\end{align}
 Suppose that on some ray
$\Re z=c \Im z$ $\,(c \in \BR, \, \Im z>0)$ the equality
\begin{align}     \label{5.20}&
\|\vp(z)-\wh \vp(z)\|=O(\E^{2\I rz}) \quad (\Im z \to \infty) \quad {\mathrm{for}}
\,\, {\mathrm{all}} \quad 0<r<r_0\leq l
\end{align}
holds. Then we have
\begin{align}     \label{5.21}&
v(x)=\wh v(x), \qquad 0<x<r_0.
\end{align}
\end{Tm}
\begin{proof}. Since, according to  \eqref{2.17}, Weyl functions are non-expansive, 
 we see that the matrix function $\E^{-2\I rz}\big(\vp(z)-\wh \vp(z)\big)$ is bounded
on  the line $\Im z =M+\ve\,$ $(\ve>0)$ and the inequality
\begin{align}     \label{5.22}&
\|\E^{-2\I rz}\big(\vp(z)-\wh \vp(z)\big)\|\leq 2 \E^{2r|z|} \quad (\Im z >M>0)
\end{align}
holds. Furthermore, formula \eqref{5.20} implies that
$\E^{-2\I rz}\big(\vp(z)-\wh \vp(z)\big)$ is bounded on the ray $\Re z=c \Im z$.
Therefore, applying the Phragmen-Lindel\"of Theorem with the angles
generated by the line $\Im z =M+\ve$ and the ray $\Re z=c \Im z$ ($\Im z \geq M+\ve$),
we get
\begin{align}     \label{5.23}&
\|\E^{-2\I rz}\big(\vp(z)-\wh \vp(z)\big)\|\leq M_1 \quad \mathrm{for} \quad \Im z \geq M+\ve.
\end{align}
 Because of formula \eqref{hea}, of its analog for $\wh \vp$,
$\wh \Phi_1$ and of the inequality \eqref{5.23}, we have
\begin{align}     \label{5.24}&
 \int_0^r \E^{2\I(x- r)z}\big(\Phi_1(x)-\wh \Phi_1(x)\big)dx   
=O(1/\Im(z)) \quad (\Im(z) \to + \infty)
\end{align}
uniformly with respect to $\Re(z)$.
Hence, the left-hand side of \eqref{5.24} is bounded in $\BC$
 and tends to zero on some rays. Thus, we obtain
\begin{align}     \label{5.25}&
\int_0^r \E^{2\I(x- r)z}\big(\Phi_1(x)-\wh \Phi_1(x)\big)dx   \equiv 0,
\quad {\mathrm{i.e.}}, \quad  \Phi_1(x)\equiv \wh \Phi_1(x) \quad
(0<x<r).
\end{align}
Since \eqref{5.25} holds for all $r<r_0$, we get $ \Phi_1(x)\equiv \wh \Phi_1(x)$ for $0<x<r_0$. In view of Theorem \ref{TmIP}, the last identity implies \eqref{5.21}.
\end{proof}

\section{Inverse problem on the semi-axis and evolution of the Weyl function
for NLS} \label{IP}
\setcounter{equation}{0}
\paragraph{1.}
In this section, we again consider Dirac system \eqref{1.1} on the semi-axis $[0,\, \infty)$ with a
bounded potential $v$: 
\begin{align}\label{IR15}
 \|v(x)\| \leq M  \quad \mathrm{for} \,\, x\in \BR_+.
 \end{align}
 According to Definition \ref{defWeyl}, the Weyl function
 of   system \eqref{1.1} such that \eqref{IR15} holds is given by the right-hand side
 of \eqref{2.3}. From Proposition \ref{PnW1} we see that the Weyl function $\vp$ is unique.
 Moreover, it follows from the proof of Corollary \ref{cyHea} that 
\begin{align}\label{IR16}
\vp(z)=2\I z\int_0^{\infty}\E^{2\I x z}\Phi(x)dx, \quad \E^{-2(M+\ve)x}\Phi(x)\in L^2_{m_2\times m_1},
 \end{align}
 and for any  $0<l<\infty$ we have
 \begin{align}\label{IR17}
 \Phi(x)\equiv \Phi_1(x) \quad (0<x<l).
 \end{align}
 Since $\Phi_1\equiv \Phi$, we get that $ \Phi_1(x)$ does not depend
on $l$ for $l>x$. Compare this with the proof of Proposition 4.1 in \cite{FKS2},
where the  fact that $E(x,t)$ (and so $\Phi_1$) does not depend on $l$
follows from the uniqueness of the factorizations of operators $S_l^{-1}$.
See also Section 3 in \cite{AGKLS2} on the uniqueness of the accelerant.
Using \eqref{IR16}, \eqref{IR17} and Theorem \ref{TmIP} we get our next result.
\begin{Tm}\label{saxis} 
The unique  Weyl function of Dirac system \eqref{1.1} satisfying \eqref{IR15} admits representation
\begin{align}\label{IR18}
\vp(z)=2\I z\int_0^{\infty}\E^{2\I x z}\Phi_1(x)dx, \quad \E^{-2\eta x}\Phi_1(x)\in L^2_{m_2\times m_1} \quad (\eta >M).
 \end{align}
The procedure to recover $\Phi_1$ and $v$ from $\vp$ is given in Theorem \ref{TmIP}.
\end{Tm}
\paragraph{2.}
The well-known   ``focusing'' nonlinear Schr\"odinger (fNLS) equation
\begin{align}  \label{e2.1} &
2v_t+\I(v_{xx}+2vv^*v)=0 \quad \Big(v_t:=\frac{\p}{\p t}v\Big),
\end{align}
where $v(x,t)$ is an $m_1 \times m_2$ matrix function,  is interesting from a mathematical point of view and plays an essential role
(in its scalar, matrix and multicomponent forms) in numerous applications (see, e.g., 
\cite{AblPr, SuSu}).
The fNLS equation \eqref{e2.1} is equivalent 
(see \cite{FT, ZS} and references therein)
to the zero curvature
equation 
\begin{align}  \label{zc} &
G_t-F_x+[G,F]=0, \quad [G,F]:=GF-FG,
\end{align}
where the $m \times m$ ($m=m_1+m_2$)
matrix functions $G(x,t,z)$ and $F(x,t,z)$ are given by the formulas
\begin{align}  \label{e2.2} &
G=\I zj+jV, \quad
F=\I\big(z^2 j-\I z jV-\big(V_x+jV^2\big)/2\big),
\end{align}
and $j$ and $V$ are defined in \eqref{1.2}.
The zero curvature representation \eqref{zc} of the integrable nonlinear equations
appeared soon after the seminal Lax pairs (see 
a historical discussion in \cite{FT} and original papers \cite{AKNS, Nov, ZM}).
Furthermore, condition \eqref{zc} is the compatibility condition for the
auxiliary linear systems
\begin{align}  \label{au} &
y_x=Gy, \quad y_t=Fy,
\end{align}
and in the case of the fNLS equation the first of  those systems is 
(in view of \eqref{e2.2}) the skew-self-adjoint
Dirac system \eqref{1.1}. The compatibility condition was studied in a more rigorous way and 
the corresponding important factorization formula for fundamental solutions was
introduced in \cite{SaL1', SaL2+}. It was proved in greater detail and under weaker
conditions in \cite{ALS11}. More specifically, we have the following proposition.
\begin{Pn} \label{TmM} \cite{ALS11}.
Let $m \times m$ matrix functions $G$ and $F$ and their derivatives $G_t$ and
$F_x$  exist on the semi-strip 
\begin{equation}      \label{e0.3}
{\mathcal D}=\{(x,\, t):\,0 \leq x <\infty, \,\,0\leq t<a\},
\end{equation}
let $G$, $G_t$ and $F$ be continuous with respect to $x$ and $t$ on $\mathcal{D}$,
 and let \eqref{zc}
hold. Then we have the equality 
\begin{equation} \label{e0.6'}
u(x,t,z)R(t,z)=R(x,t,z)u(x,0,z), \quad R(t,z):=R(0,t,z),
\end{equation}
where $u(x,t,z)$ and $R(x,t,z)$ are normalized fundamental solutions given, respectively, by:
\begin{equation} \label{e1}
u_x=Gu, \quad u(0,t,z)=I_m; \quad R_t=FR, \quad R(x,0,z)=I_m.
\end{equation}
The equality \eqref{e0.6'} means that the matrix function 
$$y=u(x,t,z)R(t,z)=R(x,t,z)u(x,0,z)$$ 
satisfies both systems \eqref{au} in 
$\cld$. Moreover, the fundamental solution $u$ admits the factorization
\begin{equation} \label{e2}
u(x,t,z)=R(x,t,z)u(x,0,z)R(t,z)^{-1}.
\end{equation}
\end{Pn}
To derive the evolution $\vp(t,z)$ of the Weyl functions  of the skew-self-adjoint
Dirac systems $y_x(x,t,z)=G(x,t,z)y(x,t,z)$, we rewrite \eqref{e2} in the form
\begin{equation} \label{e2'}
u(x,t,z)^{-1}=R(t,z)u(x,0,z)^{-1}R(x,t,z)^{-1}.
\end{equation}
Note that an additional parameter $t$ appears now in  the functions
discussed in Sections \ref{DP} and \ref{fi} (see, e.g, $u(x,t,z)$ and $\vp(t,z)$).
We partition $R(t,z)$ into blocks $R=\{R_{ij}\}_{i,j=1}^2$ (and $R_{11}$  is
an $m_1\times m_1$ block here).
  \begin{Tm}\label{evol} Let an $m_1 \times m_2$ matrix function $v(x,t)$ in $\cld$
 be  continuously differentiable together with its first derivatives
 and let  $v_{xx}$ exist. 
 Assume that $v$  satisfies the fNLS equation \eqref{e2.1}  as well as  the following inequalities:
  \begin{equation} \label{e2.8'}
  \sup_{(x,t)\in {\mathcal D}}\|v(x,t)\| \leq M, \quad \sup_{x \in \BR_+}\|v_x(x,t)\|<\infty
  \quad \mathrm{for}\,\mathrm{each} \quad 0\leq t<a.
\end{equation}
 Then the evolution $\vp(t,z)$ of the Weyl functions  of  Dirac systems $y_x=Gy$,
 where $G$ has the form \eqref{e2.2}, is given $($for $\Im(z)>M)$ by the equality
 \begin{equation} \label{e2.9}
\vp(t,z)=\big(R_{21}(t,z)+R_{22}(t,z)\vp(0,z)\big)
\big(R_{11}(t,z)+R_{12}(t,z)\vp(0,z)\big)^{-1}.
\end{equation}
  \end{Tm}
  \begin{proof}. Since $V=V^*$ and $jV =-V j$, it follows from \eqref{e2.2}
  that  $F(x,t, z)^*=-F(x,t, \ov z)$. Hence, using   \eqref{e2.2} and \eqref{e1} we get
\begin{align}&\nn
  \frac{\p}{\p t}\big(R(x,t,  z)^*jR(x,t,z)\big)=R(x,t,  z)^*(jF(x,t,z)-F(x,t, \ov z)j)R(x,t,z)
  \\ & \label{e3}=
R(x,t,  z)^*\big(\I(z^2-\ov z^2)I_m+(z+\ov z)V(x,t)+\I V_x(x,t)j\big)R(x,t,z) .
\end{align}
In view of \eqref{e2.8'} and \eqref{e3}, for each $0\leq t<a$ there are numbers $M_{1,2}(t)$
($M_1>M>0$, $M_2>0$) and the corresponding quarterplane
\begin{align}&\label{e3'}
\cld_1=\{z:\, \Im z>M_1,\,\, \Re z>M_2\} \subset \BC,
 \end{align}
such that the inequality
\begin{align}&\label{e4}
  \frac{\p}{\p t}\big(R(x,t,  z)^*jR(x,t,z)\big)<0 \end{align}
  holds in $\cld_1$.
 Inequality \eqref{e4} and the initial condition
 $R(x,0,z)=I_m$ imply $R(x,t,  z)^*jR(x,t,z)\leq j$, or, equivalently,
\begin{align} \label{e2.10}&
\big(R(x,t,  z)^*\big)^{-1}jR(x,t,z)^{-1}\geq j, \qquad z \in \cld_1.
\end{align}
Let ${\cal P}(z)$ satisfy \eqref{2.2'}  and thus satisfy \eqref{2.1}. Let $\wt\clp(x,t,z)$ (sometimes we will
also write $\wt\clp(z)$, omitting $x$ and $t$) be determined by the equality
\begin{equation} \label{e2.12}
\wt\clp(x,t,z):=R(x,t,z)^{-1}\clp(z).
\end{equation}
In view of  \eqref{e2.10}  the matrix function $\wt\clp$ satisfies \eqref{2.1} 
in $\cld_1$.

Because of the smoothness conditions on $v$, we see that $G$ and $F$ given
by \eqref{e2.2}  satisfy the requirements of Proposition \ref{TmM},
that is, \eqref{e2'} holds. Using \eqref{e2'} and \eqref{e2.12}, we see that
\begin{align} \label{e5}
u(x,t,z)^{-1}\clp(z)=
R(t,z)
\begin{bmatrix} I_{m_1} \\ \phi(x,t, z) \end{bmatrix}
\begin{bmatrix} I_{m_1} & 0 \end{bmatrix}u(x,0,z)^{-1}\wt \clp(z),
\end{align}
where
\begin{align} \label{e6}
\phi(x,t,z)=
\begin{bmatrix} 0 & I_{m_2} \end{bmatrix}u(x,0,z)^{-1}\wt \clp(z)
\left(\begin{bmatrix} I_{m_1} & 0 \end{bmatrix}u(x,0,z)^{-1}\wt \clp(z)\right)^{-1}.
\end{align}
Here, \eqref{2.7} yields inequalities
\begin{align} \label{e7}
\det \left(\begin{bmatrix} I_{m_1} & 0 \end{bmatrix}u(x,0,z)^{-1}\wt \clp(z)\right)\not=0,
\quad
\det \left(\begin{bmatrix} I_{m_1} & 0 \end{bmatrix}u(x,t,z)^{-1} \clp(z)\right)\not=0.
\end{align}
According to \eqref{e5} and \eqref{e7} we have (in the quarterplane $\cld_1$) the equality
\begin{align}\nn&
\begin{bmatrix}
0 &I_{m_2}
\end{bmatrix}u(x,t,z)^{-1}\clp(z)\Big(\begin{bmatrix}
I_{m_1} & 0
\end{bmatrix}u(x,t,z)^{-1}\clp(z)\Big)^{-1}
\\ & \label{e8}
=\big(R_{21}(t,z)+R_{22}(t,z)\phi(x,t,z)\big)
\big(R_{11}(t,z)+R_{12}(t,z)\phi(x,t,z)\big)^{-1}.
\end{align}
Taking into account Definition \ref{set}, we see that the left-hand side
of \eqref{e8} belongs $\cln(x,t,z)$ and $\phi(x,t,z) \in \cln(x,0,z)$,
where $\phi$ is defined in \eqref{e6}. Therefore, \eqref{e2.9} follows from \eqref{2.3} and \eqref{e8} (see also
\eqref{2.19}), when $x$ tends to infinity. Although, we first derived
\eqref{e2.9} only for $\cld_1$, we see that it holds everywhere in $\BC_M$ via analyticity. 
\end{proof}

\begin{Rk} \label{step} Theorems \ref{saxis} and \ref{evol} 
can be applied to recover solutions of the fNLS. Theorems on the evolution of the Weyl functions
also make up the first step in proofs of uniqueness and existence
of the solutions of nonlinear equations via the ISpT method (see, for instance, \cite{SaA7'}).
\end{Rk}


{\bf Acknowledgement.}
The work of I.Ya. Roitberg was supported by the 
German Research Foundation (DFG) under grant no. KI 760/3-1 and
the work of A.L. Sakhnovich was supported by the Austrian Science Fund (FWF) under Grant  no. Y330.



\begin{flushright} \it
B. Fritzsche,  \\
Fakult\"at f\"ur Mathematik und Informatik, \\
Mathematisches Institut, Universit\"at Leipzig, \\
Johannisgasse 26,  D-04103 Leipzig, Germany,\\
e-mail: {\tt  Bernd.Fritzsche@math.uni-leipzig.de } \\   $ $ \\

B. Kirstein, \\
Fakult\"at f\"ur Mathematik und
Informatik, \\
Mathematisches Institut, Universit\"at Leipzig,
\\ Johannisgasse 26,  D-04103 Leipzig, Germany, \\
e-mail: {\tt Bernd.Kirstein@math.uni-leipzig.de } \\  $ $ \\

I. Roitberg, \\
Fakult\"at f\"ur Mathematik und
Informatik, \\
Mathematisches Institut, Universit\"at Leipzig, \\
Johannisgasse 26,  D-04103 Leipzig, Germany, \\
e-mail: {\tt Inna.Roitberg@math.uni-leipzig.de} \\  $ $ \\

A.L. Sakhnovich, \\  Fakult\"at f\"ur Mathematik,
Universit\"at Wien,
\\
Nordbergstrasse 15, A-1090 Wien, Austria, \\
e-mail: {\tt al$_-$sakhnov@yahoo.com }
\end{flushright}

\end{document}